\documentclass[12pt]{article}
\usepackage{amssymb}
\usepackage{amsfonts}
\usepackage{mathrsfs}
\usepackage{graphicx}
\usepackage{amsmath}
\usepackage{epsfig}
\usepackage{verbatim}
\usepackage{url}

\def\smallpage{
\addtolength\textwidth{3cm} \addtolength\oddsidemargin{-1.5cm}
\addtolength\textheight{3cm} \addtolength\topmargin{-1.5cm}}
 \smallpage

\newcommand{\lab}[1]{\label{#1}}                

\newcommand{\remove}[1]{}
\newcommand\eqn[1]{(\ref{#1})}

\newcommand{\be}{\begin{equation}}
\newcommand{\ee}{\end{equation}}
\newcommand{\bea}{\begin{eqnarray}}
\newcommand{\eea}{\end{eqnarray}}
\newcommand{\bean}{\begin{eqnarray*}}
\newcommand{\eean}{\end{eqnarray*}}

\newtheorem{thm}{Theorem}[section]
\newtheorem{cor}[thm]{Corollary}

\newtheorem{lemma}[thm]{Lemma}

\newtheorem{prop}[thm]{Proposition}
\def\proof{\noindent{\bf Proof.\ }  }
\def\qed{~~\vrule height8pt width4pt depth0pt}

\def\G{{\mathcal G}}
\def\M{{\mathcal M}}

\def\hE{{\widehat{\bf E}}}
\def\EM{{\bf E}_{\mathcal M}}
\def\EH{{\bf E}_{\mathcal H}}
\def\PM{{\bf P}_{\mathcal M}}
\def\PH{{\bf P}_{\mathcal H}}
\def\po{{\bf Po}}
\def\ex{{\mathbb E}}
\def\pr{{\mathbb P}}
\def\HH{{\mathcal H}}
\def\D{{\mathcal D}}

\def\hmu{{\widehat\mu}}
\def\K{{\mathcal C}}
\def\tY{{\widetilde Y}}

\def\cD{{\overline{\mathcal D}}}
\def\bmu{{\overline\mu}}
\def\hf{{\widehat f}}

\def\eps{\epsilon}
\def\la{\lambda}

\def\ss{\smallskip}
\def\non{\nonumber}
\def\no{\noindent}

\catcode`@=11 \@addtoreset{equation}{section}

\catcode`@=12
\date{}
\title{The first $k$-regular subgraph is large}
\author{Pu Gao\\
University of Toronto\\
pu.gao@utoronto.ca}
\begin{document}
\maketitle
\begin{abstract}
Let $(G_m)_{0\le m\le \binom{n}{2}}$ be the random graph process starting from the empty graph on vertex set $[n]$ and with a random edge added in each step. Let $m_k$ denote the minimum integer such that $G_{m_k}$ contains a $k$-regular subgraph. We prove that for all sufficiently large $k$, there exist two constants $\eps_k\ge \sigma_k> 0$, with $\eps_k\to 0$ as $k\to\infty$, such that asymptotically almost surely any $k$-regular subgraph of $G_{m_k}$ has size between $(1-\eps_k)|\K_k|$ and $(1-\sigma_k)|\K_k|$, where $\K_k$ denotes the $k$-core of $G_{m_k}$.
\end{abstract}

\noindent AMS 2010 Mathematics subject classification: 05C80

\section{Introduction}
\lab{s:intro}
Let $\G_{n,m}$ denote the probability space of random graphs on $n$ vertices and $m$ edges with the uniform distribution and let $\G(n,p)$ denote the binomial model of random graphs on vertex set $[n]:=\{1,2,\ldots,n\}$, where each edge occurs independently with probability $p$.  This paper relates to the research on the existence of $k$-regular subgraphs of $\G_{n,cn/2}$ (correspondingly $\G(n,c/n)$), when $c$ is slightly greater than $c_k$, the threshold of the emergence of the $k$-core. It was shown in~\cite{L} that the property of having a $k$-regular subgraph has a sharp threshold $p^*(n,k)$ in $\G(n,p)$ for all $k\ge 3$, where $p^*(n,k)$ is of order $\Theta(1/n)$. Whether the limit of $np^*(n,k)$ exists is not known. Let $\overline {c_k}=\limsup_{n\to\infty} np^*(n,k)$ and $\underline {c_k}=\liminf_{n\to\infty} np^*(n,p)$. It is obvious that a $k$-regular subgraph exists only if the $k$-core is non-empty, whereas the converse is not necessarily true except for $k\le 2$. Thus, $\underline {c_k}\ge c_k$. In~\cite{BKV}, Bollob\'{a}s, Kim and Verstra\"{e}te proved that $\underline {c_3}$ is strictly greater than $c_3$. However, whether $\underline {c_k}$ is strictly greater than $c_k$, for general $k\ge 4$, is not clear. Bollob\'{a}s, Kim and Verstra\"{e}te conjectured that this strict inequality holds also for all $k\ge 4$. However, Pretti and Weigt~\cite{PW2} claimed the opposite: for all $k\ge 4$, $c_k$ and $\overline {c_k}$ are equal (thus $\underline {c_k}=\overline{c_k}=c_k$). Their analysis uses the (non-rigorous) cavity method and certain statistical physics techniques. (Basically they argued that $\overline {c_k}$ is the critical point where a certain entropy becomes zero and that $\overline {c_k}$ is a non-trivial solution of a certain equation, the solution of which coincides with $c_k$, for every $k\ge 4$. In the application of the cavity method they assumed some hypothesis on the {\em stability of replica symmetry}, for which they only managed to provide numerical evidence instead of a mathematical proof.)

The most recent progress on this topic are improvements of upper bounds of $\overline {c_k}$.
Pra\l at, Verstra\"{e}te and Wormald~\cite{PVW} proved that for all sufficiently large $k$, $\overline {c_k}\le c_{k+2}$ by showing that for any $c$ slightly greater than $c_{k+2}$, asymptotically almost surely (a.a.s.) the $(k+2)$-core contains a $k$-factor (or a $k$-regular subgraph that spans all but at most $k-1$ vertices of the $(k+2)$-core). More recently, Chan and Molloy~\cite{CM} improved it further to $\overline{c_k}\le c_{k+1}$, for all sufficiently large $k$. Thus, the currently best known bounds of $\underline{c_k}$ and $\overline {c_k}$ are $c_k\le \underline{c_k}\le \overline{c_k}\le c_{k+1}$, for large $k$.

In this paper, we study the size of the first $k$-regular subgraph in the graph evolution process defined as follows.
Let $G_0,G_1,\ldots,G_{\binom{n}{2}}$ be a graph process where $G_0$ is the empty graph on vertex set $[n]$ and for every $1\le i\le \binom{n}{2}$, $G_i$ is obtained by adding one edge $x$ to $G_{i-1}$, where $x$ is uniformly at random chosen from all edges not in $G_{i-1}$. Therefore, $G_0\subseteq G_1\subseteq \cdots G_{\binom{n}{2}}$ and for every $0\le m\le \binom{n}{2}$, $G_m$ is distributed exactly as $\G_{n,m}$. This random graph process was first introduced by Erd\H{o}s and R\'{e}nyi. For details, see~\cite{ER2}.
Let $m_k$ denote the minimum integer such that $G_{m_k}$ contains a $k$-regular subgraph. What is the size of a typical $k$-regular subgraph in $G_{m_k}$? Obviously, any $k$-regular subgraph of $G_{m_k}$ is a subgraph of $\K_k$, the $k$-core of $G_{m_k}$.
 For $k=3$, an observation by Pretti and Weigt~\cite{PW2} suggests that the first $3$-regular subgraph contains around $24\%$ of the vertices in $G_{m_3}$. (They claimed that they will identify the size of the first $k$-regular subgraph for general $k\ge 4$ in a following publication, but we did not find a paper on that.)
In this paper, we prove that for all sufficiently large $k$, asymptotically almost surely (a.a.s.) the first $k$-regular subgraph that appears in the random graph process misses at most $\eps_k |\K_k|$ vertices of $\K_k$, where $\eps_k\to 0$ as $k\to\infty$. On the other hand, with a simple second moment argument, we will show that a.a.s.\ this $k$-regular subgraph must miss at least $\Omega(|\K_k|)$ vertices of $\K_k$.

It follows immediately as a corollary that for any $c$ slightly greater than $c_k$, a.a.s.\ the sizes of all $k$-regular subgraphs of $\G_{n,cn/2}$ (or $\G(n,c/n)$), if there exists one, lie between $(1-\eps_k)|\K_k|$ and $(1-\sigma_k)|\K_k|$, where $\eps_k\ge \sigma_k>0$ and $\eps_k\to 0$ as $k\to\infty$. Note that this does not confirm the existence of a $k$-regular subgraph in $\G_{n,cn/2}$ and thus does not confirm or disprove that $\overline {c_k}=c_k$ for $k\ge 4$. Hence, either to rigorously prove Pretti and Weigt's prediction, or to prove the conjecture by Bollob\'{a}s, Kim and Verstra\"{e}te, we only need to restrict our future investigation to the existence of $k$-regular subgraphs with size in a narrow range $(1-\eps_k)|\K_k|$ and $(1-\sigma_k)|\K_k|$, where $\eps_k\ge \sigma_k>0$ and $\eps_k\to 0$ as $k\to\infty$.

\section{Main Results}
\lab{sec:results}

Let $n$, $M$ and $k$ be positive integers such that $M\ge kn$ is even. Let $\M(n,M,k)$ denote the probability space of random multigraphs on vertex set $[n]$ with $M/2$ edges whose end vertices are independently and uniformly at random (u.a.r.) chosen from $[n]$, conditional on that each vertex has degree at least $k$. Let $\HH(n,M,k)$ denote the probability space of $\M(n,M,k)$ restricted to simple graphs. It is well known (see~\cite{CW}) that the $k$-core of $\G_{n,m}$ (and $\G(n,p)$) is distributed as $\HH(n',2m',k)$, conditional on the number of vertices and edges in the $k$-core being $n'$ and $m'$.

Given a graph $G$ and a positive integer $k$, let $\K_k(G)$ denote the $k$-core of $G$. For a sequence of probability spaces indexed by $n$ (e.g.\ $\M(n,M,k)$ and $\HH(n,M,k)$), we say an event $A_n$ is true asymptotically almost surely (a.a.s.) if the probability that $A_n$ holds goes to $1$ as $n\to \infty$. For two functions $f(n)$ and $g(n)$ of $n$, we write $f(n)=O(g(n))$ if there is a constant $C>0$ such that $|f(n)|\le C|g(n)|$ for large $n$. We write $f(n)=o(g(n))$ if $\lim_{n\to\infty} f(n)/g(n)=0$. All unspecified asymptotics refer to $n\to\infty$. Some asymptotics refer to $k\to\infty$. In the latter case, we will always specify it.

In places where the parameter under discussion is clearly integral (e.g.\ $m=cn/2$ should be an integer), we omit the floor function if the error caused by omitting it does not affect the analysis.

We will prove two results on the sizes of the $k$-regular subgraphs of $\HH(n,M,k)$, when $M/n$ is in a certain range.

\begin{thm}\lab{t:lower}
Let $G\sim \HH(n,M,k)$, where $d=M/n=k+o(k)$ (as $k\to\infty$) and $d\ge k+\frac{1}{2}\sqrt{k\log k}$. Then, for every $\eps>0$, there exists $K>0$ such that for all $k\ge K$, a.a.s.\ there is no $k$-regular subgraph of $G$ whose size is between $\eps n$ and $(1-\eps)n$.
\end{thm}

\begin{thm}\lab{t:upper}
Let $k\ge 3$ be a fixed integer. Let $G\sim\HH(n,M,k)$, where $M=O(n)$. Then there exists $\eps>0$, such that a.a.s.\ there is no $k$-regular subgraph of $G$ whose size is more than $(1-\eps)n$.
\end{thm}

Now we consider the random graph process $(G_m)_{0\le m\le\binom{n}{2}}$ defined in Section~\ref{s:intro}, and recall that $m_k$ denotes the minimum integer such that $G_{m_k}$ contains a $k$-regular subgraph.
In order to apply Theorems~\ref{t:lower} and~\ref{t:upper} and show that the same conclusions hold also for the $k$-core of $G_{m_k}$, we need to show that the average degree of the $k$-core of $G_{m_k}$ satisfies the hypotheses of Theorems~\ref{t:lower} and Theorem~\ref{t:upper}.

For $k\ge 3$, define
\be
f_k(\mu)=\sum_{i\ge k}e^{-\mu}\frac{\mu^i}{i!}, \quad \mbox{and}\quad h_{k}(\mu)=\frac{\mu}{f_{k-1}(\mu)},\lab{fk}
\ee
and let
$$
c_{k}=\inf\{h_{k}(\mu),\mu>0\}.
$$
For any $c\ge c_{k}$, define $\mu_{c,k}$ to be the larger solution of $h_{k}(\mu)=c$. (There are two solutions for any $c>c_k$ and there is a unique solution when $c=c_k$.) In particular, let $\mu_{k}=\mu_{c_{k},k}$ for $k\ge 3$ and let $d_k=\mu_k f_{k-1}(\mu_k)/f_{k}(\mu_k)$.
Pittel, Spencer and Wormald~\cite{PSW} determined the (a.a.s.) size and density of the $k$-core for any $c>c_k$. We cite their result as follows (in a less precise form).

\begin{thm}\lab{t:coreThreshold}
Let $k \ge 3$ be fixed. Suppose that $c\le c_{k}-n^{-1/3}$ and $m=cn/2$. Then a.a.s.\ $\G_{n,m}$ has an empty $k$-core. Suppose
$c\ge c_{k}+n^{-1/3}$ and $m=cn/2$. 
Then, a.a.s.\ $\G_{n,m}$ has a non-empty $k$-core with $f_k(\mu_{c,k})n+o(n)$ vertices and $\frac{1}{2}\mu_{c,k}f_{k-1}(\mu_{c,k})n+o(n)$ edges. The same conclusions hold for $\G(n,p)$ with $p=c/n$.
\end{thm}

In the same paper, they estimated the size of the first $k$-core in the random graph process $(G_i)_{0\le i\le\binom{n}{2}}$ (see~\cite[Theorems 1 and 3]{PSW}), from which we can easily deduce the following lemma (by noting that $\mu f_{k-1}(\mu)/f_k(\mu)$ is an increasing function on $\mu>0$).

\begin{lemma} \lab{l:coreThreshold}
A.a.s.\ the number of vertices in $\K_k(G_{m_k})$ is $f_k(\mu_k)n+o(n)$ and the average degree of $\K_k(G_{m_k})$ is $d_k+o(1)$.
\end{lemma}

It was determined also in~\cite{PSW} that $c_k=k+\sqrt{k\log k}+O(\sqrt{k/\log k})$. (The error term in~\cite{PSW} is $O(\log k)$, and was corrected in~\cite{PVW}.) A more precise expression of $c_k$ was given in~\cite[Lemma 1]{PVW}, from which we can easily deduce that $c_k< 3k$ for every $k\ge 3$. In the following lemma, we estimate $d_k$ and $\mu_k$ and show that they are both close to $c_k$ when $k$ is large.

\begin{lemma}\lab{l:ck}
As $k\to\infty$,
$$
d_k=k+\sqrt{k\log k}+O(\sqrt{k/\log k}),\quad \mu_k=k+\sqrt{k\log k}+O(\sqrt{k/\log k})\ \mbox{and}\ f_{k}(\mu_k)\to 1.
$$
\end{lemma}

By Lemmas~\ref{l:coreThreshold} and~\ref{l:ck}, we see that a.a.s.\ the average degree of $\K_k(G_{m_k})$ is at least $k+\frac{1}{2}\sqrt{k\log k}$, as required by Theorem~\ref{t:lower}, whereas the condition of its average degree being at most $k+o(k)$ can be easily verified by noting that $\overline {c_k}\le c_{k+1}$ for large $k$, proved in~\cite{CM}.
It is easy to prove that a.a.s.\ for all $k\ge 3$, $\K_k(G_{m_k})$ cannot have a $k$-regular subgraph with at most $\eps_0 |\K_k|$ vertices, for some small fixed $\eps_0>0$ (by applying Lemma~\ref{l:smallS} in Section~\ref{sec:main}).  Then, using Theorems~\ref{t:lower} and~\ref{t:upper} and by Lemma~\ref{l:ck}, we can prove the following main theorem for the sizes of the $k$-regular subgraphs of $G_{m_k}$.

\begin{thm}\lab{t:main}
For all sufficiently large $k$, there exist $\eps_k\ge \sigma_k>0$ such that $\eps_k\to 0$ as $k\to\infty$ and a.a.s.\ all $k$-regular subgraphs of $G_{m_k}$ have size between $(1-\eps_k)|\K_k|$ and $(1-\sigma_k)|\K_k|$, where $\K_k=\K_k(G_{m_k})$.
\end{thm}

We also have the following result for $\G_{n,m}$ and $\G(n,p)$.



\begin{thm}\lab{t2:main}
Let $k\ge 3$. For every constant $c_k< c< 3k$ and $c=k+o(k)$ (with respect to $k\to\infty$), there exist two constants $\eps_k\ge\sigma_k>0$, with $\eps_k\to 0$ as $k\to\infty$, such that a.a.s.\ all $k$-regular subgraphs of $\G_{n,cn/2}$ (and $\G(n,c/n)$), if they exist, have size between $(1-\eps_k)|\K_k|$ and $(1-\sigma_k)|\K_k|$.
\end{thm}

We will prove Theorem~\ref{t:lower} in Section~\ref{sec:lower} and Theorem~\ref{t:upper} in Section~\ref{sec:upper}.
We provide proofs of Lemma~\ref{l:ck} and Theorems~\ref{t:main} and~\ref{t2:main} in Section~\ref{sec:main}.

\section{Proof of Theorem~\ref{t:lower}}
\lab{sec:lower}

 Let ${\bf Y}$ denote the degree sequence of $\M(n,M,k)$. It was proved in~\cite{CW} that the distribution of ${\bf Y}$ is precisely the truncated multinomial $Multi(n,M,k)$ defined as follows. Let $\D_{n,M,k}$ denote the set of vectors
$$
\left\{{\bf d}=(d_i)_{i=1}^n:\ \sum_{1\le i\le n}d_i=M,\quad d_i\ge k,\forall\ 1\le i\le n\right\}.
$$
Then for any ${\bf d}\in \D_{n,M,k}$,
\begin{equation}
\pr({\bf Y}={\bf d})=\frac{\binom{M}{d_1,\ldots,d_n}/n^M}{\sum_{{\bf d}\in\D_{n,M,k}} \binom{M}{d_1,\ldots,d_n}/n^M}=\frac{\prod_{i=1}^n 1/d_i!}{\sum_{{\bf d}\in\D_{n,M,k}} \prod_{i=1}^n 1/d_i!}.\lab{degProb}
\end{equation}

The distribution of ${\bf Y}$ can be approximated by $n$ independent copies of a truncated Poisson random variable, defined as follows. Let $d=M/n$. Choose $\la>0$ such that $\la f_{k-1}(\la)=d f_k(\la)$. It is easy to see that $\la$ exists and is unique as long as $d>k$, by noting that $\la f_{k-1}(\la)/f_k(\la)$ is an increasing function of $\la>0$. Define $Z_{\ge k}(\la)$ to be the random variable with probability function
\be
\pr(Z_{\ge k}(\la)=j)=\frac{e^{-\la}\la^j}{f_k(\la)j!}, \ \ \forall j\ge k,\lab{Zla}
\ee
where $f_k(\la)=\sum_{i\ge k}e^{-\la}\la^i/i!$ is as defined in~\eqn{fk}.

The following result on approximating ${\bf Y}$ by independent copies of $Z_{\ge k}(\la)$ can be found in~\cite[Corollary 5.3]{GW5}.
\begin{prop}\lab{p:PoissonApprox} Assume $d=M/n=O(1)$ and let $\la$ be chosen such that $\la f_{k-1}(\la)=d f_k(\la)$.
Let $A_n$ be a subset of $\D_{n,M,k}$. Let $\pr_{TP}(A_n)$ denote the probability that $(Z_1,\ldots,Z_n)\in A_n$ where $Z_i$ are independent copies of $Z_{\ge k}(\la)$ defined in~\eqn{Zla} and let $\PM(A_n)$ denote the probability that $(Y_1,\ldots, Y_n)\in A_n$ in $\M(n,M,k)$. Assume that $M-kn\to\infty$ as $n\to\infty$. Then
$$
\PM(A_n)=O(\sqrt{M})\pr_{TP}(A_n).
$$
\end{prop}

Given a degree sequence ${\bf d}\in \D_{n,M,k}$, let $g({\bf d})$ denote the number of simple graphs with degree sequence ${\bf d}$. In particular,
let $g_k(n)$ denote the number of $k$-regular graphs on $n$ vertices.
The estimation of $g({\bf d})$ was studied in a few research papers. See~\cite{M2, MW2, MW3}.
Here we cite the result by McKay~\cite{M2}. Let $d_{\max}=\max\{d_i:\ 1\le i\le n\}$. Recall that ${\bf d}\in \D_{n,M,k}$. So $\sum_{i=1}^n d_i=M$. If $d_{\max}=o(M^{1/4})$, we have the following asymptotic estimate:
\begin{equation}
g({\bf d})=\frac{M!}{2^{M/2}(M/2)!\prod_{i=1}^n d_i!}\exp\left(-\varphi({\bf d})\right), \lab{graphCount}
\end{equation}
where
\be
\varphi({\bf d})=\sum_{i=1}^n d_i(d_i-1)/2M+\left(\sum_{i=1}^n d_i(d_i-1)/2M\right)^2+O(d_{\max}^4/M).\lab{phi}
\ee
The proof of~\eqn{graphCount} uses the configuration model, first introduced by Bollob\'{a}s~\cite{B6}. Consider each vertex $i$ as a bin containing $d_i$ points. Take a uniformly random matching of all $M$ points and represent each pair in the matching as an edge in the resulting (multi)graph. Then the total number of matchings is $M!/2^{M/2}(M/2)!$. The resulting (random) graph is not necessarily simple. However, it is easy to see that every simple graph corresponds to exactly $\prod_{i=1}^n d_i!$ distinct matchings. Hence, the above estimate was obtained by proving that the probability that the resulting graph is simple is $\exp\left(-\varphi({\bf d})\right)$, which holds when $d_{\max}=o(M^{1/4})$.
Of course, if we ignore the probability of the resulting graph being simple, we have the following coarse upper bound of $g(\bf d)$
 \be
 g({\bf d})\le\frac{M!}{2^{M/2}(M/2)!\prod_{i=1}^n d_i!},\lab{roughCount}
 \ee
which holds for any degree sequence ${\bf d}$.

 Let $\PM(\cdot)$ and $\EM(\cdot)$ denote the probability and expectation in the probability space $\M(n,M,k)$ and let $\PH(\cdot)$ and $\EH(\cdot)$ denote the probability and expectation in $\HH(n,M,k)$. Let $G\sim\M(n,M,k)$. Then, for any event $A$ and any random variable $X$,
$$
\PH(A)=\frac{\PM(A \wedge (G\ simple))}{\PM(G\ simple)},\quad\EH(X)=\frac{\EM(X I\{G\ simple\})}{\PM(G\ simple)},
$$
as $\HH(n,M,k)$ is $\M(n,M,k)$ restricted to simple graphs. The following proposition is a standard method of proving a.a.s.\ properties in $\HH(n,M,k)$. Instead of proving that some property holds a.a.s.\ in $\HH(n,M,k)$ directly, it is usually easier to prove that it holds a.a.s.\ in $\M(n,M,k)$ instead. This proposition can be found in many papers (e.g.~\cite{GW5}). (In fact, it can be easily deduced from~\eqn{graphCount} and Lemma~\ref{l:badSeq} below.)
\begin{prop}\lab{p:simple}
Assume $M=O(n)$. Then $\PM(G\ simple)=\Omega(1)$, where $G\sim\M(n,M,k)$.
\end{prop}

Let $\la$ be such that $\la f_{k-1}(\la)=d f_k(\la)$, where $d=M/n$. Define
\be
\zeta=\la^2 f_{k-2}(\la)/2df_k(\la).\lab{zeta}
\ee
Define $\D_{n,M,k}^0$ to be the subset of $\D_{n,M,k}$ satisfying
\be
d_{\max}\le M^{1/4}/\log n,\quad \zeta-1/\log n\le\sum_{i=1}^n d_i(d_i-1)/2M\le \zeta+1/\log n.\lab{D0}
\ee

\begin{lemma}\lab{l:badSeq}  Assume $M=O(n)$. Then
$\PM({\bf Y}\notin \D_{n,M,k}^0)=o(1)$.
\end{lemma}
\proof Let $Z_1,\ldots,Z_n$ be $n$ independent copies of $Z_{\ge k}(\la)$, where $\la f_{k-1}(\la)=d f_k(\la)$ and $d=M/n$. Then $\ex Z_i=\la^2f_{k-2}(\la)/f_{k}(\la)$. First we bound the probability
$$
\pr\Big(\max_{1\le i\le n} Z_i\ge M^{1/4}/\log n\Big).
$$
By the definition of $Z_{\ge k}(\la)$ in~\eqn{Zla} and putting $t=M^{1/4}/\log n$, we immediately have
\be
\pr\Big(\max_{1\le i\le n} Z_i\ge t\Big)\le \sum_{i=1}^n \pr(Z_i\ge t)=O(n e^{-\la}\la^{t}/f_{k}(\la)t!)=O(n(e\la /t)^t)=O(e^{-t}).\lab{maxDeg}
\ee
Next, we bound the probability that
$$
\left|\sum_{i=1}^n Z_i(Z_i-1)/2M-\zeta\right|>1/\log n.
$$
Let $W_i=Z_iI\{Z_i\le\log n\}$. Then for every $1\le i\le n$,
\begin{eqnarray*}
\pr(Z_i\neq W_i)&=&\pr(Z_i>\log n)=O(e^{-\la}\la^{\log n}/(\log n)!)<\exp(-\la+\log n \log(e\la/\log n))\\
&<&\exp\left(-\frac{1}{2}\log n\log\log n\right)=n^{-\log\log n/2}.
\end{eqnarray*}
We have
$$
\ex W_i(W_i-1)=\la^2 f_{k-1}(\la)/f_k(\la)=2d\zeta.
$$
By Azuma-Hoeffding's inequality~\cite[Theorem 2.25]{JLR} and noting that $M=dn$ and $|W_i|\le \log n$ for all $i$, we have
\begin{eqnarray*}
&&\pr\left(\left|\sum_{i=1}^n W_i(W_i-1)/2M-\zeta\right|>\frac{1}{\log n}\right)=\pr\left(\left|\sum_{i=1}^n \Big(W_i(W_i-1)-\ex W_i(W_i-1)\Big)\right|>\frac{2M}{\log n}\right)\\
&&\qquad\le2\exp\left(-\frac{2(2M/\log n)^2}{(\log n)^2 n}\right)=2\exp\left(-\frac{8d^2n}{(\log n)^4 }\right).
\end{eqnarray*}
Hence,
\begin{eqnarray*}
&&\pr\left(\left|\sum_{i=1}^n Z_i(Z_i-1)/2M-\zeta\right|>\frac{1}{\log n}\right)\le\pr\left(\left|\sum_{i=1}^n W_i(W_i-1)/2M-\zeta\right|>\frac{1}{\log n}\right)+n\pr(Z_1\neq W_1)\\
&&\qquad\le 2\exp\left(-\frac{8d^2n}{(\log n)^4 }\right)+n^{1-\log\log n/2}=o(1/M).
\end{eqnarray*}
Hence, by Proposition~\ref{p:PoissonApprox}, $\PM({\bf Y}\notin \D_{n,M,k}^0)=o(1)$. \qed
\ss

\remove{
\begin{lemma} There exists a constant $\alpha>0$, such that for all $k$,
a.a.s.\ $\M(n,M,k)$ contains no $k$-regular subgraph with size at most $\alpha n$.
\end{lemma}
\proof For any $s$, let $X_s$ denote the number of subgraphs of $\M(n,M,k)$ with $s$ vertices and at least $t=ks/2$ edges. There are $\binom{n}{s}$ ways to choose a subset of vertices of size $s$. If $s\ge\log^2 n$, given a set $S$ of $s$ vertices, the probability that S contains more than $1.1 ds$ points is at most $\beta^{ds}$ for some $0<\beta<1$ by Lemma~\ref{l:degConcentration}. Conditional on $S$ containing at most $1.1ds$ points. The probability that $S$ contains at least $t$ pairs with both ends inside $S$ is at most
$$
\binom{1.1ds}{2t}\frac{(2t)!}{t!2^t}\prod_{i=0}^{t-1}\frac{1}{D-1-2i}\le \left(\frac{e(1.1ds)^2}{2t(D-2t)}\right)^t
$$
}

In the rest of this section, let $\eps>0$ be a small but fixed constant and
let $B_{\eps}$ denote the event that there is a $k$-regular subgraph whose size is between $\eps n$ and $(1-\eps)n$. Given a set $S$, let $A_S$ denote the event that there is a $k$-regular subgraph on $S$. Then, by Proposition~\ref{p:simple} and Lemma~\ref{l:badSeq} and by symmetry,
\bea
\PH(B_{\eps})&\le& \frac{\PM(B_{\eps}\wedge (G\ simple) \wedge {\bf Y}\in \D_{n,M,k}^0)+\PM({\bf Y}\notin \D_{n,M,k}^0)}{\PM(G\ simple)}\non\\
&=&O(\PM(B_{\eps}\wedge (G\ simple) \wedge {\bf Y}\in \D_{n,M,k}^0))+o(1)\non\\
&=&O\left(\sum_{\eps n\le s\le (1-\eps)n}\sum_{S\subseteq [n]:\ |S|=s}\PM(A_{S}\wedge (G\ simple) \wedge {\bf Y}\in \D_{n,M,k}^0)\right)+o(1)\non\\
&=&O\left(\sum_{\eps n\le s\le (1-\eps)n}\binom{n}{s}\PM(A_{[s]}\wedge (G\ simple) \wedge {\bf Y}\in \D_{n,M,k}^0)\right)+o(1).\lab{B}
\eea

Let $X_{[s]}$ denote the number of $k$-regular subgraphs on $S=[s]$. For an arbitrary subset $\D$ of $\D_{n,M,k}^0$, let $\cD=\D_{n,M,k}^0\setminus \D$. Then, by the Markov inequality,
\bea
&&\PM(A_{[s]}\wedge (G\ simple) \wedge {\bf Y}\in \D_{n,M,k}^0)=\PM(A_{[s]}\wedge (G\ simple) \wedge {\bf Y}\in \D)+O(\PM(\cD))\non\\
&&\quad\le \EM (X_{[s]}I\{(G\ simple)\wedge {\bf Y}\in \D\})+O(\PM(\cD))\non\\
&&\quad=\sum_{{\bf d}\in \D}\EM(X_{[s]} \mid (G\ simple)\wedge {\bf Y}={\bf d})\PM(G\ simple\mid {\bf Y}={\bf d} )\PM({\bf Y}={\bf d})+O(\PM(\cD))\non\\
&&\quad\le\sum_{{\bf d}\in \D}\EH(X_{[s]} \mid  {\bf Y}={\bf d})\PM({\bf Y}={\bf d})+O(\PM(\cD))\non\\
&&\quad\le\sum_{{\bf d}\in \D} g_k(s)\frac{g({\bf d'})}{g({\bf d})}\PM({\bf Y}={\bf d})+O(\PM(\cD)),\lab{A}
\eea
where ${\bf d}'=(d'_i)_{i=1}^n$ is defined by
$d_i'=d_iI\{i>s\}+(d_i-k)I\{i\le s\}$ for all $1\le i\le n$.
Note that the last inequality above holds because the number of $k$-regular graphs on $S$ is $g_k(s)$, whereas the probability that a given $k$-regular graph on $S$ is a $k$-regular subgraph of $\HH(n,M,k)$ on $S$, conditional on the degree sequence of $\HH(n,M,k)$ being ${\bf d}$, is at most $g({\bf d'})/g({\bf d})$.

In what follows, we will choose appropriate $\D$ so that $\PM(\cD)$ is sufficiently small and we will upper bound
\be
\sum_{{\bf d}\in \D} g_k(s)\frac{g({\bf d'})}{g({\bf d})}\PM({\bf Y}={\bf d}).\lab{sum}
\ee

Let $M_s\ge ks$ and let $\M(n,M,M_{s},k)$ denote the probability space of $\M(n,M,k)$ conditioned to $\sum_{i\le s}Y_i=M_{s}$.
 Then the distribution of the degree sequence $Y_1\ldots, Y_s$, a subsequence of ${\bf Y}$, is precisely the truncated multinomial $Multi(s,M_s,k)$. 
Let
$$
\D_{n,M,M_{s},k}=\left\{{\bf d}\in\D^0_{n,M,k}:\ \sum_{i\le s} d_i=M_{s}\right\}.
$$
Let $\sigma>0$ be fixed. Define
$$
\D_{\sigma}=\underset{ks\le M_s\le(1+\sigma)ds}{\cup}\D_{n,M,M_{s},k}.
$$

We will estimate~\eqn{sum} and $\PM(\cD)$ with $\D=\D_{\sigma}$.  First we upper bound $\PM(\cD_{\sigma})$.

\begin{lemma}\lab{l:mu-d}  Assume $d\ge k+\frac{1}{2}\sqrt{k\log k}$ and $d=k+o(k)$ (as $k\to\infty$).
Take $\mu$ that satisfies $\mu f_{k-1}(\mu)=d f_k(\mu)$. Then $k\le\mu=k+o(k)$ for all large $k$ and $\mu/d\to 1$, as $k\to\infty$.
\end{lemma}
\proof Let $\hmu=k$. Then $f_k(\hmu)\ge 1/2$. Hence, for all large $k$,
$$
\frac{\hmu f_{k-1}(\hmu)}{f_k(\hmu)}=\hmu\left(1+\frac{e^{-\hmu}\hmu^{k-1}}{(k-1)!f_k(\hmu)}\right)\le k+C\sqrt{k}\le d,
$$
for some positive constant $C$, by Stirling's approximation. It is easy to check that $x f_{k-1}(x)/f_k(x)$ is an increasing function of $x$ and so we
have that $\mu>\hmu=k$. We also observe that $\mu=d f_k(\mu)/f_{k-1}(\mu)\le d=k+o(k)$.  So $\mu/d\to 1$ as $k\to\infty$.

For a set $S$ of vertices, let $deg(S)$ denote the sum of degrees of vertices in $S$.
The following lemma is a standard concentration result on the degree sum of a set of vertices that is not too small. See~\cite[Corollary 5.4]{GW5} for a detailed proof.

\begin{lemma}\lab{l:degConcentration} Let $G\sim\M(n,M,k)$ where $M=O(n)$ and $d=M/n\ge k$. Take $\mu$ that satisfies $\mu f_{k-1}(\mu)=d f_k(\mu)$. Assume $\mu/d\to 1$ as $k\to\infty$. Then
for every $\sigma>0$, there exist $K>0$ and $0<\alpha<1$ such that for all $k>K$, and for any $S\subset [n]$ with $|S|\ge\log^2 n$,
$$
\PM(|deg(S)-d|S||>\sigma d|S|)=\alpha^{d|S|}.
$$
\end{lemma}

By Lemmas~\ref{l:mu-d} and~\ref{l:degConcentration}, for any fixed $\sigma>0$, we have
\be
\PM(\cD_{\sigma})\le \beta^{d s},\quad \mbox{for some}\ 0<\beta<1.\lab{PcD}
\ee
 Next, we estimate~\eqn{sum}. Note that for every $\sigma>0$, we have $1+\sigma\ge k/d$, as $d\ge k$. Thus, $\D_{\sigma}$ is non-empty. We will prove the following lemma.
\begin{lemma}\lab{l:sum} Let $G\sim\M(n,M,k)$ where $M=O(n)$ and $d=M/n\ge k+\frac{1}{2}\sqrt{k\log k}$. For every $\sigma>0$ and $\eta>1$,  
define
\be
f(\eta,\sigma)=\frac{\left(1-1/\eta\right)^{\eta/2}}{\sqrt{\eta-1}}\left(\frac{e^{-\mu}}{f_k(\mu)}\right)^{1/k}\mu^{(1+\sigma)d/k } \left(\frac{e}{(1+\sigma)d-k}\right)^{(1+\sigma)d/k -1} ,\lab{f} 
\ee
where $\mu$ is the root of
$$
\mu f_{k-1}(\mu)=d f_{k}(\mu).
$$
Then there is an absolute constant $C>0$ (independent of $k$) such that for every $\sigma>0$, and for every $\eps n\le s\le (1-\eps)n$,
$$
\sum_{{\bf d}\in \D_{\sigma}} g_k(s)\frac{g({\bf d'})}{g({\bf d})}\PM({\bf Y}={\bf d})\le 2ds M \Big(\big(1+C(\sigma+d/k-1)\big) f(M/ks,\sigma)\Big)^{ks}.
$$
\end{lemma}
We leave Lemma~\ref{l:sum} to be proved in Section~\ref{s:sum}. In this section, we complete the proof of Theorem~\ref{t:lower} by applying Lemma~\ref{l:sum}.

Let $\delta=s/n$. By~\eqn{A} and Lemma~\ref{l:sum}, for any fixed $\sigma>0$,
$$
\PM\Big(A_{[s]}\wedge (G\ simple) \wedge {\bf Y}\in \D_{n,M,k}^0\Big)\le 2dsM\big((1+C(\sigma+d/k-1))f(d/k\delta,1+\sigma)\big)^{ks}+O(\PM(\cD_{\sigma})).
$$
Since $\binom{n}{s}\le (\delta^{-\delta}(1-\delta)^{-1+\delta})^n$ eventually for every $\eps n\le s\le (1-\eps)n$, by~\eqn{B}, we have
\bea
\PH(B_{\eps})&=&O\left(\sum_{\eps n\le s\le (1-\eps)n} 2dsM\Big(\delta^{-1/k}(1-\delta)^{-\frac{1-\delta}{k\delta}}(1+C(\sigma+d/k-1)))f(d/k\delta,\sigma)\Big)^{ks}\right)\nonumber\\
&+&O\left(\sum_{\eps n\le s\le (1-\eps)n}\binom{n}{s}\PM(\cD_{\sigma})\right)+o(1).\lab{B2}
\eea
\begin{lemma}\lab{l:PcD} Assume $d\ge k+\frac{1}{2}\sqrt{k\log k}$. For every $\eps>0$ and every $\sigma>0$, there exists a $K>0$ such that for all $k>K$,
$$
\sum_{\eps n\le s\le (1-\eps)n}\binom{n}{s}\PM(\cD_{\sigma})=o(1).
$$
\end{lemma}

\proof By~\eqn{PcD}, there exists $0<\beta<1$ depending only on $\sigma$ that
$$
\sum_{\eps n\le s\le (1-\eps)n}\binom{n}{s}\PM(\cD_{\sigma})\le\sum_{\eps n\le s\le (1-\eps)n}\binom{n}{s}\beta^{ds} \le\sum_{s\ge \eps n} \left(\frac{en\beta^d}{s}\right)^s\le\sum_{s\ge \eps n} \left(\frac{e\beta^d}{\eps}\right)^s.
$$
Choose $K$ sufficiently large so that $e\beta^d/\eps<1$. Then the above summation is $o(1)$.\qed\ss

Recall the definition of $f$ in~\eqn{f}, where $\mu$ is the root of $\mu f_{k-1}(\mu)=d f_k(\mu)$.
\begin{lemma} \lab{l:less1} Assume $d\ge k+\frac{1}{2}\sqrt{k\log k}$ and $d=k+o(k)$ (as $k\to\infty$). Then, for every $\eps>0$, there exist sufficiently large $K>0$ and sufficiently small constants $\sigma>0$ and $\hat\eps>0$ such that for all $k>K$ and for all $\eps\le \delta\le 1-\eps$,
$\delta^{-1/k}(1-\delta)^{-\frac{1-\delta}{k\delta}}(1+C(\sigma+d/k-1)))f(\eta,\sigma)<1-\hat\eps$, where $\eta=d/k\delta$.
\end{lemma}

\proof By Lemma~\ref{l:mu-d}, $\mu\ge k$ and so $f_k(\mu)\ge 1/2$. Then, clearly, as $k\to\infty$,
$$
\delta^{-1/k}(1-\delta)^{-\frac{1-\delta}{k\delta}}\to 1, \quad \frac{e^{1/k}}{k!^{1/k}}=(1+o(1)) e/k, \quad  \frac{\mu}{k}\to 1, \quad f_k(\mu)^{1/k}\to 1.
$$
Define $g(\eta,\sigma)=\delta^{-1/k}(1-\delta)^{-\frac{1-\delta}{k\delta}}(1+C(\sigma+d/k-1)))f(\eta,\sigma)$. (Note that $\delta=d/k\eta$ depends only on $\eta$.) Then, as $k\to\infty$ and $\sigma\to 0$,
$$
g(\eta,\sigma)\to (\eta-1)^{\eta/2-1/2} \eta^{-\eta/2} e^{-\mu/k+(1+\sigma)d/k}\left(\frac{\mu}{(1+\sigma) d-k}\right)^{(1+\sigma) d/k-1},
$$
since $\sigma+d/k-1\to 0$ as $k\to\infty$ and $\sigma\to 0$.
Next, we show that
$$
e^{-\mu/k+(1+\sigma) d/k}\left(\frac{\mu}{(1+\sigma) d-k}\right)^{(1+\sigma) d/k-1}\to 1,
$$
as $k\to\infty$ and $\sigma\to 0$.
Let $x_{\mu}=\mu-k$ and $x_d=d-k$.  Then, $x_{\mu}=o(k)$ and $x_{d}=o(k)$ by the assumption that $d=k+o(k)$ and Lemma~\ref{l:mu-d}. So,
\begin{eqnarray*}
&&e^{-\mu/k+(1+\sigma) d/k}\left(\frac{\mu}{(1+\sigma) d-k}\right)^{(1+\sigma) d/k-1} =e^{(x_d-x_{\mu})/k+\sigma(1+x_d/k)}\left(\frac{k+x_{\mu}}{x_d+\sigma(k+x_d)}\right)^{x_d/k+\sigma(1+x_d/k)}\\
&&\quad=\exp\Big((1+o(1))\sigma\Big) \left(\frac{k+x_{\mu}}{x_d+\sigma(k+x_d)}\right)^{x_d/k+\sigma(1+x_d/k)}\\
&&\quad=\exp\Big((1+o(1))\sigma\Big)\exp\left(\left(\frac{x_d}{k}+\sigma(1+x_d/k)\right)\log\left(\frac{1+x_{\mu}/k}{x_d/k+\sigma(1+x_d/k)}\right)\right),
\end{eqnarray*}
where the asymptotics above refers to $k\to\infty$.
Since $x_{\mu}/k\to 0$ and $x_d/k+\sigma(1+x_d/k)\to 0$, as $k\to\infty$ and $\sigma\to 0$, the above can be arbitrarily small by choosing sufficiently large $k$ and sufficiently small $\sigma$.
Hence, $g(\eta,\sigma)\to (\eta-1)^{\eta/2-1/2} \eta^{-\eta/2}$. It is easy to prove that $(\eta-1)^{\eta/2-1/2} \eta^{-\eta/2}$ is a strictly decreasing function on $\eta>1$ with limit $1$ as $\eta$ approaches to $1$ from above. Since $d>k$ and $\delta\le 1-\eps$, we have $\eta=d/k\delta\ge 1/(1-\eps)$. Thus, there exists $\eps'>0$, such that $(\eta-1)^{\eta/2-1/2} \eta^{-\eta/2}<1-\eps'$ and so $g(\eta,\sigma)<1-\hat\eps$, for some small constant $\hat\eps>0$, by choosing sufficiently small $\sigma$ and sufficiently large $k$ (so that $x_d/k$ and $x_{\mu}/k$ are sufficiently small).
Therefore, for any $\eps>0$, there exist sufficiently large $K>0$ and sufficiently small $\sigma>0$, such that for all $k>K$, $\delta^{-1/k}(1-\delta)^{-\frac{1-\delta}{k\delta}}(1+C(\sigma+d/k-1)))f(\eta,\sigma)<1-\hat\eps$ for all $\eps\le\delta\le 1-\eps$.\qed
\ss

Now we complete the proof of Theorem~\ref{t:lower}. By Lemma~\ref{l:less1}, for every $\eps>0$, there exist constants $\sigma>0$ and $K>0$, such that for all $k>K$,
$$
\sum_{\eps n\le s\le (1-\eps)n} 2dsMn\Big(\delta^{-1/k}(1-\delta)^{-\frac{1-\delta}{k\delta}}(1+C(\sigma+d/k-1)))f(\eta,\sigma)\Big)^{ks}< \sum_{\eps n\le s\le (1-\eps)n} \beta^{ks}=o(1),
$$
for some constant $0<\beta<1$. Let $\sigma>0$ be chosen so that the above holds. Then by Lemma~\ref{l:PcD}, provided $k$ is sufficiently large, $\sum_{\eps n\le s\le (1-\eps)n}\binom{n}{s}\PM(\cD_{\sigma})=o(1)$. By~\eqn{B2}, $\PH(B_{\eps})=o(1)$ and so a.a.s.\ there is no $k$-regular subgraph in $\HH(n,M,k)$ with size between $\eps n$ and $(1-\eps)n$.

\subsection{Proof of Lemma~\ref{l:sum}}
\lab{s:sum}
Recall that for any ${\bf d}\in \D_{n,M,M_{s},k}$, ${\bf d}'$ is defined as $d_i'=d_iI\{i>s\}+(d_i-k)I\{i\le s\}$. By~\eqn{roughCount} and~\eqn{graphCount},
\begin{eqnarray*}
g_k(s)&\le&\frac{(ks)!}{2^{ks/2}(ks/2)! k!^{s}},\\
g({\bf d})&=&\frac{M!}{2^{M/2}(M/2)!\prod_{i=1}^n d_i!}\exp(-\varphi({\bf d})+o(d_{\max}^4/M)), \\
g({\bf d}')&\le&\frac{(M-ks)!}{2^{(M-ks)/2}((M-ks)/2)!\prod_{i=1}^n (d_i-k)!},
\end{eqnarray*}
where $\varphi({\bf d})$ is defined in~\eqn{phi}. By the definition of $\D^0_{n,M,k}$ in~\eqn{D0}, and the definition of $\zeta$ in~\eqn{zeta},
for any ${\bf d}\in \D_{n,M,M_{s},k}\subseteq \D^0_{n,M,k}$, we have
$$
\varphi({\bf d})=\zeta+\zeta^2+o(1).
$$
Then for some constant $C>0$,

\begin{eqnarray}
&&\hspace{-.8cm}\sum_{{\bf d}\in \D_{\sigma}} g_k(s)\frac{g({\bf d'})}{g({\bf d})}\PM({\bf Y}={\bf d})=\sum_{M_s=ks}^{(1+\sigma)d s}\sum_{{\bf d}\in \D_{n,M,M_s,k}} g_k(s)\frac{g({\bf d'})}{g({\bf d})}\PM({\bf Y}={\bf d})\nonumber\\
&&\le \frac{(ks)!\exp(\zeta+\zeta^2+o(1))}{2^{ks/2}(ks/2)! k!^{s}}\frac{2^{M/2}(M/2)!(M-ks)!}{M!2^{(M-ks)/2}((M-ks)/2)!}\sum_{M_s=ks}^{(1+\sigma)d s}\sum_{{\bf d}\in \D_{n,M,M_s,k}}\prod_{1\le i\le s}[d_i]_k\PM({\bf Y}={\bf d})\nonumber\\
&&\le C\left(\frac{ks}{M-ks}\right)^{ks/2}\frac{1}{k!^{s}}\left(1-\frac{ks}{M}\right)^{M/2}\sum_{M_s=ks}^{(1+\sigma)d s}\sum_{{\bf d}\in \D_{s,M_s,k}}\prod_{1\le i\le s}[d_i]_k\PM({\bf \tY}={\bf d}),\lab{summation}
\end{eqnarray}
where ${\bf \tY}=(Y_1,\ldots, Y_s)$ is a subsequence of ${\bf Y}$. Since the distribution of ${\bf \tY}$ is $Multi(s,M_s,k)$, by~\eqn{degProb},
\begin{eqnarray*}
\sum_{{\bf d}\in \D_{s,M_s,k}}\prod_{1\le i\le s}[d_i]_k\PM({\bf \tY}={\bf d})&=&\sum_{{\bf d}\in \D_{s,M_{s},k}} \frac{\prod_{i\le s}[d_i]_k\prod_{i=1}^{s} 1/d_i!}{\sum_{{\bf d}\in\D_{s,M_{s},k}} \prod_{i=1}^{s} 1/d_i!}\\
&=&\frac{1}{\sum_{{\bf d}\in\D_{s,M_{s},k}} \prod_{i=1}^{s} 1/d_i!}\sum_{{\bf d}\in \D_{s,M_{s},k}} \prod_{i=1}^{s} \frac{1}{(d_i-k)!}.
\end{eqnarray*}

\remove{
Let $\hE(\cdot)$ denote $\EM(\cdot\mid (G\ simple)\wedge\sum_{i\le s} Y_i=M_{s})$. Then
\begin{eqnarray*}
\hE(X_S)&=&g_k(s)\sum_{{\bf d}\in \D_{n,S,M,M_{s},k}}\frac{g({\bf d'})}{g({\bf d})}\pr({\bf Y}={\bf d})\\
&=&\frac{(ks)!}{2^{ks/2}(ks/2)! k!^{s}}\frac{2^{M/2}(M/2)!(M-ks)!}{M!2^{(M-ks)/2}((M-ks)/2)!}\sum_{{\bf d}\in \D_{s,M_{s},k}} \frac{\prod_{i\le s}[d_i]_k\prod_{i=1}^{s} 1/d_i!}{\sum_{{\bf d}\in\D_{s,M_{s},k}} \prod_{i=1}^{s} 1/d_i!}\\
&=&\frac{(ks)!(M/2)!(M-ks)!}{(ks/2)! k!^{s}M!((M-ks)/2)!}\frac{1}{\sum_{{\bf d}\in\D_{s,M_{s},k}} \prod_{i=1}^{s} 1/d_i!}\sum_{{\bf d}\in \D_{s,M_{s},k}} \prod_{i=1}^{s} \frac{1}{(d_i-k)!}.
\end{eqnarray*}
}

Obviously, for every $M_s\ge ks$, $\prod_{i=1}^{s} 1/d_i!$ for ${\bf d} \in\D_{s,M_s,k}$ is minimized when ${\bf d}=(k,k,\ldots,k,k+\phi)$, where $\phi=M_s-ks$. 
Thus, we immediately have the following lemma.
\begin{lemma}\lab{l:sum0} For every $M_s\ge ks$,
\begin{eqnarray*}
&&\sum_{{\bf d}\in\D_{s,M_{s},k}} \prod_{i=1}^{s} \frac{1}{d_i!}\ge k!^{-s+1}(k+\phi)!^{-1}=k!^{-s}\frac{k!}{(k+\phi)!}\ge \frac{k!^{-s}}{(k+\phi)^{\phi}},
\end{eqnarray*}
where $\phi=M_s-ks$.
\end{lemma}

When $M_s-ks\to\infty$, we have a fairly precise estimate of $\sum_{{\bf d}\in\D_{s,M_{s},k}} \prod_{i=1}^{s} 1/d_i!$, stated in the following lemma.
\begin{lemma} \lab{l:sum1} Suppose $M_s-ks\to\infty$ as $n\to\infty$. Then
$$
\sum_{{\bf d}\in\D_{s,M_{s},k}} \prod_{i=1}^{s} 1/d_i!=\Theta\Big(1/\sqrt{M_{s}}\Big)e^{\mu s}\mu^{-M_{s}}f_k(\mu)^{s},
$$
where $\mu$ satisfies
$$
\frac{\mu f_{k-1}(\mu)}{f_k(\mu)}=\frac{M_{s}}{s}.
$$
\end{lemma}
\proof
Let $Z_i$, $i\le s$ be independent copies of the truncated Poisson variables $Z_{\ge k}(\mu)$, defined in~\eqn{Zla}.
Then,
$$
\pr\Big(\sum_{i\le s}Z_i=M_{s}\Big)=\sum_{{\bf d}\in \D_{s,M_{s},k}}\prod_{i\le s}\frac{e^{-\mu} \mu^{d_i}}{f_k(\mu)d_i!}=\frac{e^{-\mu s}\mu^{M_{s}}}{f_k(\mu)^{s}}\sum_{{\bf d}\in \D_{s,M_{s},k}}\prod_{i\le s}\frac{1}{d_i!}.
$$
By the definition of $\mu$, $\pr\big(\sum_{i\le s}Z_i=M_{s}\big)=\Theta\big(1/\sqrt{M_{s}}\big)$ by~\cite[Theorem 4(a)]{PW}. Hence,
$$
\sum_{{\bf d}\in \D_{s,M_{s},k}}\prod_{i\le s}\frac{1}{d_i!}=\Theta\Big(1/\sqrt{M_{s}}\Big)e^{\mu s}\mu^{-M_{s}}f_k(\mu)^{s}.\qed
$$


In the next lemma we deduce an upper bound of $\sum_{{\bf d}\in \D_{s,M_{s},k}} \prod_{i=1}^{s} 1/(d_i-k)!$.

\begin{lemma}\lab{l:sum2} For every $M_s\ge ks$,
$$
\sum_{{\bf d}\in \D_{s,M_{s},k}} \prod_{i=1}^{s} \frac{1}{(d_i-k)!}\le\left(\frac{e}{M_{s}/s-k}\right)^{M_{s}-ks}.
$$
\end{lemma}

\proof First we observe that
$$
\sum_{{\bf d}\in \D_{s,M_{s},k}} \prod_{i=1}^{s} \frac{1}{(d_i-k)!}=\sum_{{\bf d}\in \D_{s,M_{s}-ks,0}} \prod_{i=1}^{s} \frac{1}{d_i!}.
$$

If $M_s=ks$, the inequality holds trivially as both sides equal to $1$ (by defining $0^0=1$). Suppose $M_s>ks$.
Let $Z_i$, $i\le s$ be independent copies of the Poisson variable $\po(\la)$, where $\la =M_{s}/s-k$. Then
$$
\pr\Big(\sum Z_i=M_{s}-ks\Big)=\sum_{{\bf d}\in \D_{s,M_{s}-ks,0}}\prod_{i\le s}\frac{e^{-\la} \la^{d_i}}{d_i!}=e^{-M_{s}+ks}\la^{M_{s}-ks}\sum_{{\bf d}\in \D_{s,M_{s}-ks,0}} \prod_{i\le s} \frac{1}{d_i!}.
$$
Since $\pr\big(\sum Z_i=M_{s}-ks\big)\le 1$, we have
$$
\sum_{{\bf d}\in \D_{s,M_{s}-ks,0}} \prod_{i\le s} \frac{1}{d_i!}\le \frac{e^{M_{s}-ks}}{\la^{M_{s}-ks}}.\qed
$$

By Lemmas~\ref{l:sum0},~\ref{l:sum1} and~\ref{l:sum2} and~\eqn{summation} and using Stirling's formula, we have
\begin{eqnarray}
\sum_{{\bf d}\in \D_{\sigma}} g_k(s)\frac{g({\bf d'})}{g({\bf d})}\PM({\bf Y}={\bf d})&\le& \phi+ M\left(\frac{ks}{M-ks}\right)^{ks/2}\frac{1}{k!^{s}}\left(1-\frac{ks}{M}\right)^{M/2}\non\\
&&\times \sum_{M_s=ks+\log n}^{(1+\sigma)d s}\left(\frac{e^{-\bmu}}{f_k(\bmu)}\right)^{s}\bmu^{M_{s}} \left(\frac{e}{M_{s}/s-k}\right)^{M_{s}-ks}, \lab{Bprob}
\end{eqnarray}
where $\bmu$ satisfies
$$
\frac{\bmu f_{k-1}(\bmu)}{f_k(\bmu)}=\frac{M_{s}}{s},
$$
and
\bea
\hspace{-.7cm}\phi&=&
(\log n)  M \left(\frac{ks}{M-ks}\right)^{ks/2}\frac{1}{k!^s}\left(1-\frac{ks}{M}\right)^{M/2}k!^s(k+\log n)^{\log n}\frac{e^{M_{s}-ks}}{(M_{s}/s-k)^{M_{s}-ks}}.\lab{phi}
\eea
Note that $\phi$ corresponds to the contribution to the summation~\eqn{Bprob} from $ks\le M_s<ks+\log n$.

\remove{
Let $\eps>0$ be fixed. For any fixed $\sigma>0$, let $A_{\sigma}$ denote the event that for every $S$ with $|S|\ge \eps n$,  $(1-\sigma) \delta M\le\sum_{i\in S} Y_i\le (1+\sigma)\delta M$, where $\delta=s/n$. Let $B_S$ denote the event that there is a $k$-regular subgraph on $S$, and let $B(s)$ denote the event that there is a $k$-regular subgraph on $s$ vertices. Then
\begin{eqnarray*}
\pr(B(s) \wedge A_{\sigma})&\le& \binom{n}{s}\pr(B_S\wedge A_{\sigma})=\binom{n}{s}\sum_{(1-\sigma) \delta M\le M_s\le (1+\sigma) \delta M}\pr(B_S\mid \sum_{i\in S}Y_i=M_s)\pr(\sum_{i\in S}Y_i=M_s)\\
&\le&\binom{n}{s}\max\left\{\hE(X_S),\quad (1-\sigma) \delta M\le M_s\le (1+\sigma) \delta M\right\}.
\end{eqnarray*}

Next, we estimate $\max\left\{\hE(X_S),\quad (1-\sigma) \delta M\le M_s\le (1+\sigma) \delta M\right\}$.
}

We will prove that the summand in~\eqn{Bprob} maximizes at $M_s=(1+\sigma)ds$. Before that, we first prove two technical lemmas.

\begin{lemma}\lab{l:rho}
Let $\rho> k$ and let $\mu=\mu(\rho)$ be the root of
$
\mu f_{k-1}(\mu)/f_k(\mu)=\rho
$.
Then $\mu> \rho-k$ for all $\rho> k$.
\end{lemma}
\proof Let $\hmu=\rho-k>0$.
Since $f_k(\hmu)> e^{-\hmu} \hmu^k/k!$, we have
\bean
\frac{\hmu f_{k-1}(\hmu)}{f_k(\hmu)}&=&\hmu+\frac{\hmu e^{-\hmu}\hmu^{k-1}}{(k-1)! f_k(\hmu)}=\hmu+\frac{ k e^{-\hmu}\hmu^{k}}{k! f_k(\hmu)}< \hmu+k=\rho.
\eean
Since $\mu f_{k-1}(\mu)/f_k(\mu)$ is an increasing function of $\mu$, we have $\mu> \hmu=\rho-k$. This completes the proof of the lemma.\qed\ss

\begin{lemma}\lab{l:muDer}
Let $\mu=\mu(\rho)$ be the root of
$
\mu f_{k-1}(\mu)/f_k(\mu)=\rho
$,
for $\rho>k$.
Then
$$
\mu^{\rho} \frac{e^{-\mu}}{f_k(\mu)} \left(\frac{e}{\rho-k}\right)^{\rho-k}
$$
is an increasing function on $\rho>k$.
\end{lemma}
\proof
Define $g(x)=xf_{k-1}(x)/f_k(x)$. Then $g(x)$ is an increasing function on $x>0$ and $g(x)\to k$ as $x\to 0$. Let $h_1(\rho)=\rho\log \mu-\mu-\log f_k(\mu)$. Then
\[
\mu^{\rho} \frac{e^{-\mu}}{f_k(\mu)}=\exp(h_1(\rho)).
\]
Taking the derivative of $h_1$ with respect to $\rho$, we have
$$
h_1'(\rho)=\log\mu+\frac{\rho}{\mu}\mu'(\rho)-\mu'(\rho)-\frac{f_k'(\mu)}{f_k(\mu)}\mu'(\rho)=\log\mu+\mu'(\rho)\left(\frac{\rho}{\mu}-1-\frac{f_k'(\mu)}{f_k(\mu)}\right).
$$
Since $f_k'(\mu)=f_{k-1}(\mu)-f_k(\mu)$ for all $k\ge 1$, we have
$$
h_1'(\rho)=\log\mu+\mu'(\rho)\left(\frac{\rho}{\mu}-\frac{f_{k-1}(\mu)}{f_k(\mu)}\right).
$$
Since $\rho=g(\mu)$, we have $\rho/\mu=f_{k-1}(\mu)/f_k(\mu)$, which implies that $h_1'(\rho)=\log\mu$.

On the other hand,
$$
\left(\frac{e}{\rho-k}\right)^{\rho-k}=\exp\Big((\rho-k)\big(1-\log (\rho-k)\big)\Big).
$$
Let $h_2(\rho)=(\rho-k)(1-\log (\rho-k))$. Then the derivative of $h_2$ is
$
-\log (\rho-k)
$.

Thus,
\[
\mu^{\rho} \frac{e^{-\mu}}{f_k(\mu)} \left(\frac{e}{\rho-k}\right)^{\rho-k}=\exp(h_1(\rho)+h_2(\rho)),
\]
and
$$
h_1'(\rho)+h_2'(\rho)=\log\mu-\log(\rho-k)=\log\left(\frac{\mu}{\rho-k}\right).
$$
By Lemma~\ref{l:rho}, the above is greater than $0$ for all $\rho>k$. Hence, $
\mu^{\rho} \frac{e^{-\mu}}{f_k(\mu)} \left(\frac{e}{\rho-k}\right)^{\rho-k}
$ is an increasing function of $\rho$ on $\rho>k$. \qed
Now, we upper bound the summation in~\eqn{Bprob}. By Lemma~\ref{l:muDer},
$$
\left(\frac{e^{-\bmu}}{f_k(\bmu)}\right)^{s}\bmu^{M_{s}} \frac{e^{M_{s}-ks}}{(M_{s}/s-k)^{M_{s}-ks}}=\left(\mu^{\rho} \frac{e^{-\bmu}}{f_k(\bmu)} \left(\frac{e}{\rho-k}\right)^{\rho-k}\right)^s
$$
 in~\eqn{Bprob} (taking $\rho=M_s/s$) is maximized at $M_s=(1+\sigma)ds$.
 Let $\eta=M/ks$. 
 Define
$$
\hf(\eta,\sigma)=\frac{\left(1-1/\eta\right)^{\eta/2}}{\sqrt{\eta-1}}\frac{1}{k!^{1/k}}\left(\frac{e^{-\hmu}}{f_k(\hmu)}\right)^{1/k}\hmu^{(1+\sigma) d/k } \left(\frac{e}{(1+\sigma)d-k}\right)^{(1+\sigma)d/k-1},
$$
where $\hmu$ is the root of
$$
\hmu f_{k-1}(\hmu)/f_{k}(\hmu)=(1+\sigma)d .
$$
Compare $\hf$ with $f$ in~\eqn{f}. They are almost the same except that $\mu$ and $\hmu$ are defined differently.
By~\eqn{Bprob},
$$
\sum_{{\bf d}\in \D_{\sigma}} g_k(s)\frac{g({\bf d'})}{g({\bf d})}\PM({\bf Y}={\bf d})\le \phi+(1+\sigma)d s M \hf(\eta,\sigma)^{ks},
$$
where $\phi$ is defined in~\eqn{phi}.
\remove{

 Define
\bean
g(\eta)&=&\max\left\{k!^{1/k}(k+\log n)^{\log n/ks},\hmu^{\alpha \delta \eta }\left(\frac{e^{-\hmu}}{f_k(\hmu)}\right)^{1/k}\right\},\\
\hf(\eta)&=&\frac{1}{\sqrt{\eta-1}}\frac{e^{1/k}}{k!^{1/k}}\left(1-1/\eta\right)^{\eta/2}g(\eta) ,
\eean
where $\hmu$ is the root of
$$
\hmu f_{k-1}(\bmu)/f_{k}(\hmu)=(1+\sigma)d .
$$
Compare $\hf$ with $f$ in~\eqn{f}. They are almost the same except that $\mu$ and $\hmu$ are slightly different and there is an extra term $k!^{1/k}(k+\log n)^{\log n/ks}$ in the definition of $\hf$.
By~\eqn{Bprob} and Corollary~\ref{c:muDer} ( by noting that $\phi\le M \hf(\eta)^{ks}$),
\be
\sum_{{\bf d}\in \D_{\sigma}} g_k(s)\frac{g({\bf d'})}{g({\bf d})}\PM({\bf Y}={\bf d})\le 2ds M \hf(\eta)^{ks}.\lab{sum2}
\ee
}
Now we complete the proof of Lemma~\ref{l:sum}. Comparing the terms in~\eqn{phi} with terms in~\eqn{Bprob}, It is straightforward to verify that eventually
$$
k!^{1/k}(k+\log n)^{\log n/ks}\frac{e^{M_{s}/ks-1}}{(M_{s}/s-k)^{M_{s}/ks-1}}\le \left(\frac{e^{-\hmu}}{f_k(\hmu)}\right)^{1/k}\hmu^{(1+\sigma) d/k } \frac{e^{(1+\sigma)d/k -1}}{\big((1+\sigma)d-k\big)^{(1+\sigma)d/k-1}}
$$
for every $ks\le M_s\le ks+\log n$, since $s\to\infty$ as $n\to\infty$. Thus,
$$
\sum_{{\bf d}\in \D_{\sigma}} g_k(s)\frac{g({\bf d'})}{g({\bf d})}\PM({\bf Y}={\bf d})\le 2ds M \hf(\eta,\sigma)^{ks},
$$
since $(1+\sigma)ds+\log n\le 2ds$ eventually.

By the definition of $\hf$, we only need to show that
there is $C>0$ such that $\hf(\eta,\sigma)\le(1+C(\sigma+d/k-1)))f(\eta,\sigma)$ for every $\sigma>0$ and $d=k+o(k)$. %
%
%
By Lemma~\ref{l:mu-d} and the assumption that $d\ge k+\frac{1}{2}\sqrt{k\log k}$, we have $\mu\ge k$. By the definition of $\hmu$, we have $\mu\le \hmu=(1+\sigma)d f_k(\hmu)/f_{k-1}(\hmu)\le (1+\sigma)d$. Thus,
$$
1\le \frac{\hmu}{\mu}\le (1+\sigma)\frac{d}{k} \le 1+C'(\sigma+d/k-1),
$$
for some constant $C'>0$ that does not depend on $k$ or $\sigma$, since $d=k+o(k)$.
Thus,
\bean
\hmu^{(1+\sigma)d/k}&=&(1+O(\sigma+d/k-1))\mu^{(1+\sigma)d/k}\\
\frac{e^{-\hmu/k}}{f_k(\hmu)^{1/k}}&=&(1+O(\sigma+d/k-1))\frac{e^{-\mu/k}}{f_k(\mu)^{1/k}},
\eean
where
the constants involved in the asymptotics above are independent of $\sigma$ and $k$.
This completes the proof of Lemma~\ref{l:sum}.
\section{Proof of Theorem~\ref{t:upper}}
\lab{sec:upper}
\begin{lemma} \lab{l:propB}
Let $M=O(n)$ and $M\ge kn$. For every $k\ge 3$, a.a.s.\ there exists $\Omega(n)$ vertices in $\HH(n,M,k)$ whose degrees are at least $k+1$ and whose neighbours all have degree exactly $k$.
\end{lemma}
\proof As $M=O(n)$, a.a.s.\ the numbers of vertices with degree $k$ and $k+1$ are both $\Omega(n)$ by applying Proposition~\ref{p:PoissonApprox}. I.e., there are positive constants $C_1<C_2$ and $D_1<D_2$ such that a.a.s.
\be
 C_1 n\le N_k\le C_2 n, \quad D_1 n\le N_{k+1}\le D_2 n,\lab{size}
 \ee
 where $N_i$ denote the number of vertices with degree $i\in\{k,k+1\}$. Recall the definition of $\D_{n,M,k}^0$ in~\eqn{D0}. Let $\D^1_{n,M,k}$ denote the subset of $\D_{n,M,k}^0$ such that $C_1 n\le \sum_{i:d_i=k} 1\le C_2 n$ and $D_1 n\le \sum_{i:d_i=k+1} 1\le D_2 n$. Let ${\bf d}\in \D_{n,M,k}^1$. Conditional on the degree sequence of $\HH(n,M,k)$ being ${\bf d}$, $\HH(n,M,k)$ can be generated by the configuration model, as described below~\eqn{phi}. Let $\HH_{{\bf d}}$ denote the random (multi)graph generated by the configuration model. We prove that the claim in this lemma holds in $\HH_{\bf d}$ for every ${\bf d}\in \D_{n,M,k}^1$. Then the lemma follows by~\eqn{size} and Lemma~\ref{l:badSeq} and the fact that the probability that $\HH_{{\bf d}}$ is simple is $\Omega(1)$ (See~\cite{M2}). Now consider $\HH_{{\bf d}}$ where ${\bf d}\in \D^1_{n,M,k}$. Let $n_k=\sum_{i:d_i=k} 1$ and $n_{k+1}=\sum_{i:d_i=k+1} 1$. Let $X$ denote the number of vertices with degree $k+1$, whose neighbours all have degree $k$. There are $n_{k+1}$ ways to choose a vertex $v$ with degree $k+1$ and $\binom{n_k}{k+1}$ ways to choose $k+1$ vertices as neighbours of $v$. For each of these $k+1$ chosen vertices, there are $k$ ways to choose a point inside the vertex. The number of ways to match those $k+1$ points to the $k+1$ points inside $v$ is $(k+1)!$. Then,
$$
\ex X=n_{k+1}\binom{n_k}{k+1}k^{k+1}(k+1)!\prod_{i=0}^{k}\frac{1}{M-1-2i}=\Omega(n),
$$
where $\prod_{i=0}^{k}\frac{1}{M-1-2i}$ is the probability that a given set of $k+1$ pairs occurs in the random matching over $M$ points.
Next, we compute $\ex X(X-1)$. I.e.\ we compute the expected number of ordered pairs of distinct vertices, both with degree $k+1$, and both with all neighbours having degree $k$. We first count such pairs $(v_1,v_2)$ such that $v_1$ is not adjacent to $v_2$ and $N(v_1)\cap N(v_2)=\emptyset$, where $N(v)$ denotes the set of neighbours of $v$. There are $n_{k+1}(n_{k+1}-1)$ ways to choose $(v_1,v_2)$ and $\binom{n_k}{k+1}\binom{n_{k}-k-1}{k+1}$ ways to choose the neighbours of $v_1$ and $v_2$. Hence, the expected number of such pairs is
\be
n_{k+1}(n_{k+1}-1)\binom{n_k}{k+1}\binom{n_{k}-k-1}{k+1}k^{2(k+1)}(k+1)!^2\prod_{i=0}^{2k+1}\frac{1}{M-1-2i}=(1+o(1))(\ex X)^2.\lab{0}
\ee
Next, we bound the number of such pairs $(v_1,v_2)$ such that $v_1$ is adjacent to $v_2$. There are $n_{k+1}(n_{k+1}-1)$ ways to choose $(v_1,v_2)$ and at most $\binom{n_k}{k}\binom{n_{k}}{k}$ ways to choose the neighbours of $v_1$ and $v_2$. Hence, the expected number of such pairs is at most
\be
O\left(n_{k+1}^2n_k^{2k}\prod_{i=0}^{2k}\frac{1}{M-1-2i}\right)=O(n).\lab{1}
\ee
Last, we bound the number of such pairs $(v_1,v_2)$ such that $v_1$ is not adjacent to $v_2$ and $|N(v_1)\cap N(v_2)|=j$, $j\ge 1$. There are $\binom{n_k}{j}\binom{n_{k}-j}{k+1-j}\binom{n_{k}-k-1}{k+1-j}\le n_k^{2k+2-j}$ to choose the neighbours of $v_1$ and $v_2$ and hence the expected number of such pairs is
\be
O\left(n_{k+1}^2n_k^{2k+2-j}\prod_{i=0}^{2k+1}\frac{1}{M-1-2i}\right)=O(n^{2-j})=O(n).\lab{2}
\ee
Combining~\eqn{0},~\eqn{1} and~\eqn{2}, we have $\ex X(X-1)=(1+o(1))(\ex X)^2$.
It follows then that the variance of $X$ is
$$
\ex X(X-1)+\ex X-(\ex X)^2=o((\ex X)^2).
$$
By Chebyshev's inequality, it follows that a.a.s.\ $X=\Omega(n)$ and the claim of the lemma follows thereby.\qed
\ss

\no {\bf Proof of Theorem~\ref{t:upper}.\ } Let $G\in \HH(n,M,k)$. A vertex in $G$ has property $B$ if it has degree at least $k+1$ and all its neighbours have degree exactly $k$. If $v$ has property $B$, then any $k$-regular subgraph $H$ of $G$ must miss at least one of its neighbour. By Lemma~\ref{l:propB}, a.a.s.\ there are $\Omega(n)$ vertices having property $B$. On the other hand, each vertex with degree $k$ can be adjacent to at most $k$ vertices that have property $B$. Hence, any $k$-regular subgraph of $G$ must miss $\Omega(n)$ vertices of $G$.\qed

\remove{
Then
$$
d=\frac{\mu f_{k-1}(\mu)}{f_k(\mu)},\quad c=\frac{\mu}{f_{k-1}(\mu)}.
$$

If $c=c_k$, then
$$
\frac{ke^{-\mu} \mu^k/k!f_k(\mu)}{d}=1/(k-1).
$$
It follows then that
$$
\frac{e^{-\mu}\mu^{k-1}}{(k-2)!}=f_{k-1}(\mu).
$$
It is also known that $c=k+\sqrt{k\log k}+o(\sqrt{k})$.

Numerical calculation gives
$$
3.383634281<\mu_4<3.383634283,\quad 5.075451424<d_4<5.075451426.
$$
}

\section{Proofs of Theorems~\ref{t:main} and~\ref{t2:main}}
\lab{sec:main}
Before approaching Theorems~\ref{t:main} and~\ref{t2:main}, we first give the proof of Lemma~\ref{l:ck}.
\ss

\no {\bf Proof of Lemma~\ref{l:ck}.\ }By the definition of $\mu_{k}$, $h_k'(\mu_{k})=0$. Since
$$
h'_{k}(x)=\frac{f_{k-1}(x)-xf'_{k-1}(x)}{f_{k-1}(x)^{2}},
$$
and $f'_k(x)=f_{k-1}(x)-f_k(x)$ for all $k\ge 1$, we have
$$
f_{k-1}(\mu_{k})=\mu_{k}(f_{k-2}(\mu_{k})-f_{k-1}(\mu_{k})),
$$
i.e.,
$$
\frac{f_{k-2}(\mu_{k})}{f_{k-1}(\mu_{k})}=\frac{1+\mu_{k}}{\mu_{k}}=1+\frac{1}{\mu_{k}}.
$$
On the other hand,
$$
f_{k-2}(\mu_{k})=f_{k-1}(\mu_{k})+e^{-\mu_{k}}\frac{\mu_{k}^{k-2}}{(k-2)!}.
$$
It follows immediately that
$$
\frac{\mu_k e^{-\mu_{k}}\mu_{k}^{k-2}}{(k-2)!f_{k-1}(\mu_{k})}=1.
$$
Multiply both sides by $1/(k-1)$, we get
\be
\frac{e^{-\mu_{k}}\mu_{k}^{k-1}}{(k-1)!f_{k-1}(\mu_{k})}=\frac{1}{k-1}.\lab{eq:muk}
\ee
The left hand side is a decreasing function of $\mu_k$ on $[k,+\infty)$. Moreover,
taking $\hmu=k$, we have $f_{k-1}(\hmu)\ge 1/2$ and so,
$$
\frac{e^{-\hmu}\hmu^{k-1}}{(k-1)!f_{k-1}(\hmu)}\le 2\frac{e^{-k}k^k}{k!}\le C/\sqrt{k},
$$
for some constant $C>0$ by Stirling's approximation. Since $C/\sqrt{k}$ is larger than $1/(k-1)$ for all large $k$, we have
$\mu_k\ge \hmu=k$ for all large $k$.
So,
 $f_{k-1}(\mu_k)\ge 1/2$.

Let $x=\mu_k-k$ ($x\ge 0$ as $\mu_k\ge k$). Then using Stirling's formula, we have
\begin{eqnarray*}
f_{k-1}(\mu_k)&=&\mu_k\frac{e^{-\mu_{k}}\mu_{k}^{k-2}}{(k-2)!}\sim \frac{(k+x)e^{-x-2}}{\sqrt{2\pi k}}\left(1+\frac{x+2}{k-2}\right)^{k-2}\\
&=&\frac{(k+x)e^{-x-2}}{\sqrt{2\pi k}}\exp\left(x+2-\frac{(x+2)^2}{2(k-2)}+O(x^3/k^2)\right)\\
&=&\frac{\sqrt{k}(1+O(x/k))}{\sqrt{2\pi }}\exp\left(-\frac{(x+2)^2}{2(k-2)}+O(x^3/k^2)\right).
\end{eqnarray*}
Thus,
$$
-\frac{(x+2)^2}{2(k-2)}+O(x^3/k^2)=-\frac{1}{2}\log k+\log(\sqrt{2\pi}f_{k-1}(\mu_k)(1+O(x/k))).
$$
It follows immediately that $\mu_k=k+\sqrt{k\log k}+O(\sqrt{k/\log k})$.
Then, by~\eqn{eq:muk} and the fact that $f_{k-1}(\mu_k),f_{k}(\mu_k)=\Omega(1)$,
 $$
 d_k=\mu_k\left(1+\frac{e^{-\mu_k}\mu_k^{k-1}}{(k-1)!f_k(\mu_k)}\right)=\mu_k(1+O(1/k))=\mu_k+O(1)=k+\sqrt{k\log k}+O(\sqrt{k/\log k}).
 $$
Since  $\mu_k=k+\sqrt{k\log k}+O(\sqrt{k/\log k})$, it is very easy to verify that $f_{k}(\mu_k)\to 1$ as $k\to\infty$. \qed
\ss

\begin{lemma}\lab{l:smallS} Let $k\ge 3$ and $\eps_0=1/30e^5$. Assume $p\le 3k/n$. Then, a.a.s.\ all subgraphs of $\G(n,p)$ with size $s\le \eps_0 n$ have less than $ks/2$ edges.
\end{lemma}

\proof Given $s$, let $X_s$ denote the number of sets $S$ with $|S|=s$ that contain at least $ks/2$ edges. There are $\binom{n}{s}$ ways to choose a set of $s$ vertices and given $S$ with $|S|=s$, the probability that there are at least $ks/2$ edges inside $S$ is at most
$$
\binom{\binom{s}{2}}{ks/2}p^{ks/2}\le \left(\frac{es^2p}{ks}\right)^{ks/2}\le\left(3es/n\right)^{3s/2}.
$$
Hence, using $\binom{n}{s}\le (en/s)^s$,
$$
\ex X_s\le\binom{n}{s}\left(3es/n\right)^{3s/2}\le (3e^2)^s\left(3es/n\right)^{s/2}=(27e^5s/n)^{s/2}.
$$
By the choice of $\eps_0$, the expected number of sets $S$ with $|S|\le \eps_0 n$ that contain at least $k|S|/2$ edges is at most
$$
\sum_{s\le\log n} \left(O(s/n)\right)^{s/2} + \sum_{\log n<s\le \eps_0 n} 3e^2\left(\frac{3}{4}\right)^{s/2}=o(1).
$$
The claim follows by the first moment method.\qed\ss

Recall that $m_k$ is the minimum integer such that $G_{m_k}$ contains a $k$-regular subgraph in the graph evolution process.
We will use the following theorem to upper bound $m_k$.
\begin{thm}[Chan and Molloy~\cite{CM}]\lab{t:k+1}
For all sufficiently large $k$, and for every $c_{k+1}<c<c_{k+1}+2\sqrt{k\log k}$, a.a.s.\ the $(k+1)$-core of $\G(n,c/n)$ contains a $k$-factor or a $k$-regular subgraph that expands all but one vertex.
\end{thm}

\no {\bf Proof of Theorem~\ref{t:main}.\ } By Theorems~\ref{t:coreThreshold} and~\ref{t:k+1} and the monotonicity that if $G_i$ contains a $k$-regular subgraph, then so does $G_j$ for all $j\ge i$, for all sufficiently large $k$, a.a.s.\ $\widehat m_1\le m_k\le \widehat m_2$, where $\widehat m_1=(c_k-n^{-1/3})n/2$ and $\widehat m_2=(c_{k+1}+1)n/2$. Thus, a.a.s.\ $G_{\widehat m_1}\subseteq G_{m_k}\subseteq G_{\widehat m_2}$. Let $\K_k=\K_k(G_{m_k})$. Then $\K_k$ is non-empty. By Lemma~\ref{l:coreThreshold}, the average degree of $\K_k$ is at least $d_k+o(1)$ (asymptotics referring to $n\to\infty$), which is greater than $k+\frac{1}{2}\sqrt{k\log k}$ by Lemma~\ref{l:ck}. Since a.a.s.\ $G_{m_k}\subseteq G_{\widehat m_2}\sim \G_{n,\widehat m_2}$, the average degree of $\K_k$ is a.a.s.\ $k+o(k)$ by Theorem~\ref{t:coreThreshold} (by noting that $\mu f_{k-1}(\mu)/f_k(\mu)$ is an increasing function on $\mu>0$).
Then, by Theorem~\ref{t:lower}, there are $(\eps_k)_{k\ge 3}$ with $\eps_k\to 0$ as $k\to\infty$, such that a.a.s.\ there is no $k$-regular subgraph of $\K_k$ with size between $\eps_k|\K_k|$ and $(1-\eps_k)|\K_k|$.
Since a.a.s.\ the number of edges in $\G(n,3k/n)$ is at least $\widehat m_2$ by the definition of $\widehat m_2$, we can couple $G_{\widehat m_2}$ and $\G(n,3k/n)$ so that a.a.s.\ $G_{\widehat m_2}\subseteq\G(n,3k/n)$. For details of the coupling, we refer readers to~\cite{JLR}. So, a.a.s.\ $G_{m_k}\subseteq \G_{n,\widehat m_2}\subseteq \G(n,3k/n)$. By Lemma~\ref{l:smallS}, a.a.s.\ there is no $k$-regular subgraph of $\G(n,3k/n)$ (and hence of $G_{m_k}$) with size at most $\eps_0|\K_k|\le \eps_0 n$, where $\eps_0=1/4e$. Since $\eps_k\to 0$ as $k\to\infty$, we have $\eps_0\ge \eps_k$ for all sufficiently large $k$.  Then, for all large $k$, there is no $k$-regular subgraph with size at most $(1-\eps_k)|\K_k|$. By Theorem~\ref{t:upper}, there is $\sigma_k>0$, with $\sigma_k\to 0$ as $k\to\infty$, such that a.a.s.\ there is no $k$-regular subgraph of $\K_k$ with size greater than $(1-\sigma_k)|\K|$. This completes the proof for Theorem~\ref{t:main}.
\qed \ss

\no {\bf Proof of Theorem~\ref{t2:main}.\ } Let $G\sim \G(n,c/n)$ or $G\sim \G_{n,cn/2}$.
Consider constant $c_k< c< 3k$ and $c=k+o(k)$.
 We can couple $G$ and $\G(n,3k/n)$ so that $G\subseteq \G(n,3k/n)$ since $c<3k$. Then,
by Lemma~\ref{l:smallS}, a.a.s.\ there is no $k$-regular subgraph of $G$ with size at most $\eps_0|\K_k|\le \eps_0 n$, where $\eps_0=1/30e^5$.
 Let $\mu$ be the larger solution of $\mu/f_{k-1}(\mu)=c$. Then $\mu>\mu_k$ and $\mu=k+o(k)$ and $f_{k-1}(\mu)\to 1$ as $k\to\infty$. Then, by Theorem~\ref{t:lower} and Lemma~\ref{l:smallS}, there are $(\eps_k)_{k\ge 3}$ with $\eps_k\to 0$ as $k\to\infty$, such that a.a.s.\ there is no $k$-regular subgraph of $\K_k$ with size between $\eps_k|\K_k|$ and $(1-\eps_k)|\K_k|$ (by setting $\eps_k=1-\eps_0$ for small values of $k$ that are not considered in Theorem~\ref{t:lower}).
Since $\eps_k\to 0$ as $k\to\infty$, we have $\eps_0\ge \eps_k$ for all but finitely many $k$. For each $k$ such that $\eps_0<\eps_k$, redefine $\eps_k=1-\eps_0$.  Then, for all $k\ge 3$, there is no $k$-regular subgraph with size at most $(1-\eps_k)|\K_k|$. By Theorem~\ref{t:upper}, there is $\sigma_k>0$ with $\sigma_k\to 0$ as $k\to\infty$, such that a.a.s.\ there is no $k$-regular subgraph of $\K_k$ with size greater than $(1-\sigma_k)|\K|$. This completes the proof for Theorem~\ref{t2:main}.  \qed


\begin{thebibliography}{99}
\bibitem{B6} B. Bollob\'{a}s,
Random graphs.
Second edition. {\em Cambridge Studies in Advanced Mathematics}, 73. Cambridge University Press, Cambridge, 2001. xviii+498 pp.

\bibitem{BKV}
B. Bollob\'{a}s, J. Kim and J. Verstra\"{e}te,
Regular subgraphs of random graphs,
{\em Random Structures Algorithms} 29 (2006), no. 1, 1-–13.


\bibitem{CM} S. Chan and M. Molloy, $(k+1)$-cores have $k$-factors,
{\em Combin. Probab. Comput.} 21 (2012), no. 6, 882-–896.

\bibitem{CW} J. Cain and N. Wormald, Encores on cores, {\em Electronic Journal of Combinatorics}, 13, RP 81, 2006.
\bibitem{ER2} P. Erd\H{o}s and A. R\'{e}nyi,
On the evolution of random graphs,
{\em Bull. Inst. Internat. Statist.} 38 (1961) 343--347.

\bibitem{GW5} P. Gao and N. Wormald, Orientability thresholds for random hypergraphs, preprint available at \url{http://arxiv.org/abs/1009.5489}.

\bibitem{JLR} S. Janson and T. {\L}uczak and A. Ruci{\'n}ski, Random graphs.
{\em Wiley-Interscience Series in Discrete Mathematics and Optimization}. Wiley-Interscience, New York, 2000. xii+333 pp.

\bibitem{L} S. Letzter, The property of having a k -regular subgraph has a sharp threshold, {\em Random Sturctures Algorithms} 42 (2013), pages 509--519.

\bibitem{M2}   B.D. McKay, Asymptotics for symmetric 0-1
matrices with prescribed row sums, {\em Ars Combinatoria} 19A (1985),
15--25.
\bibitem{MW2} B.D. McKay and N.C. Wormald, Asymptotic enumeration by degree
sequence of graphs with degrees $o(\sqrt{n})$, {\em Combinatorica} {\bf
11} (1991), 369--382.
\bibitem{MW3} B.D. McKay and N.C. Wormald, Asymptotic enumeration by degree sequence of graphs of high degree, {\em European journal of combinatorics}, {\bf 11}, 1990, 565--580.

\bibitem{PSW}B. Pittel and J. Spencer and N. Wormald, Sudden emergence of a giant $k$-core in a random graph, {\em J. Combin. Theory Ser. B},67, no.1, 1996, 111--151.

\bibitem{PVW} P. Pra\l at,  J. Verstra\"{e}te and N. Wormald,
On the threshold for $k$-regular subgraphs of random graphs,
{\em Combinatorica} 31 (2011), no. 5, 565-–581.

\bibitem{PW} B. Pittel and N. Wormald,
Asymptotic enumeration of sparse graphs with a minimum degree constraint,
{\em J. Combin. Theory Ser. A} 101 (2003), no. 2, 249-–263.


\bibitem{PW2}
M. Pretti and M. Weigt,
Sudden emergence of $q$-regular subgraphs in random graphs,
{\em Europhys. Lett.} 75 (2006), no. 1, 8–-14.



\end{thebibliography}
\end{document}